\documentclass{article}[10pt]

\usepackage[linesnumbered, vlined, ruled]{algorithm2e}
\usepackage{soul}
\usepackage{latexsym}
\usepackage{amssymb}
\usepackage{amsmath,amsfonts}
\usepackage{graphicx}
\usepackage{algorithm2e}
\usepackage{algorithmic}
\usepackage{todonotes}
\usepackage{amsthm}
\usepackage{geometry}
\usepackage{algorithm2e}
\usepackage{multirow}
\usepackage{hyperref}
\usepackage{subcaption}
\usepackage{float}
\restylefloat{table}

\numberwithin{equation}{subsection}

\newtheorem{theorem}{Theorem}[section]
\newtheorem{prop}[theorem]{Proposition}
\newtheorem{lemma}{Lemma}[section]
\newtheorem{corollary}{Corollary}[section]
\newtheorem{example}{Example}[subsection]
\newtheorem{definition}{Definition}[section]

\newtheorem{conjecture}[theorem]{Conjecture}

\newtheorem{question}[theorem]{Question}
\def\ex{\begin{example}}
\def\eex{\end{example}}
\def\exx{\end{example}}
\def\T{\begin{theorem}}
\def\TT{\end{theorem}}
\def\D{\begin{definition}}
\def\DD{\end{definition}}
\def\l{\begin{lemma}}
\def\ll{\end{lemma}}
\def\c{\begin{corollary}}
\def\cc{\end{corollary}}
\def\cj{\begin{conjecture}}
\def\cjj{\end{conjecture}}
\def\e{\begin{equation}}
\def\ee{\end{equation}}
\def\p{\begin{prop}}
\def\pp{\end{prop}}
\def\q{\begin{question}}
\def\qq{\end{question}}

\def \R {\mathbb{R}}

\def \inv {^{-1}}
\def \mb {\mathbf}
\DeclareMathOperator*{\argmin}{arg\,min}

\begin{document}

\baselineskip 15pt

\title{\bf Alternating Minimization for Computed Tomography with Unknown Geometry Parameters.}
\author{ Phuong Mai Huynh Pham\\ Emory University\\ phuongmai.huynhpham@gmail.com\\ \\ Manuel Santana\\ Utah State University\\ manuelarturosantana@gmail.com\\ 
\\ Advisors:\\\\ Ana Castillo\\Proximity Learning \\ana.t.castillo.rivas@gmail.com \\ \\ James Nagy\\Emory University \\ jnagy@emory.edu}
\date{}

\maketitle

\begin{abstract}
Due to the COVID-19 pandemic, there is an increasing demand for portable CT machines worldwide in order to diagnose patients in a variety of settings. This has led to a need for CT image reconstruction algorithms that can produce high quality images in the case when multiple types of geometry parameters have been perturbed. In this paper we present an alternating minimization algorithm to address this issue, where one step minimizes a regularized linear least squares problem, and the other step minimizes a bounded non-linear least squares problem. Additionally, we survey existing methods to accelerate convergence of the algorithm and discuss implementation details. Finally, numerical experiments are conducted to illustrate the effectiveness of the algorithm.
\end{abstract}

\section {Introduction}
In medical imaging, computed tomography (CT) techniques are becoming increasingly popular for their ability to produce high quality images of the human body. These help doctors diagnose several types of cancers and recently handle COVID-19 cases. 
A CT scanner is a device that is composed of a scanning gantry, X-ray generator, computer system, console panel and a physician’s viewing console. The scanning gantry is the part that produces and detects X-rays. In a typical CT scan, a patient will lay on a bed that will move through the gantry. An X-ray tube rotates around the patient and projects X-ray beams through the human body at different angles. These X-ray measurements are then processed on a computer using mathematical algorithms to create tomographic (cross-sectional) images of the tissues inside the body.  Limitations arise when using CT scanners for these medical procedures since these devices require extensive care, such as regular maintenance. Additionally, transporting these to remote locations is not an easy task. Point-of-care imaging addresses these challenges by allowing radiologists to add portable CT scanners to their departments to increase patient satisfaction and improve medical outcomes. However, the parameters that are associated with the geometry of these devices cannot always be precisely calibrated in point-of-care situations. These parameters may change over time, when the scanner is adjusted during the image acquisition process or when transported to a new location.
\begin{figure}[H]
    \centering
    \includegraphics[]{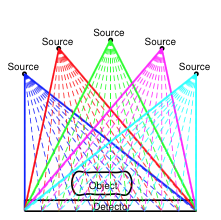}
    \caption{Object-Source-Detector}
    \label{fig:my_label}
\end{figure}

Figure 1 illustrates how an X-ray source rotates around an object during a typical CT acquisition process. The location of the source is measured by a view angle, i.e, location of the source on a circle around the object, and the distance from the source to the center of the object. If the location of the X-ray source has been perturbed, the reconstructed image will be of very poor quality. The case where only the view angles have been perturbed is an active area of research \cite{angles2} \cite{angles1}. In \cite{william} Ding devises the numerical scheme we study, and uses it to investigate when the angles and the distance from the source to the object are unknown. The purpose of this paper is to study reconstruction methods in the case when both view angles and distance from the source to the object have been perturbed. Ding performed one small experiment of this type, and we build upon his work by testing larger images, larger perturbations, more acceleration techniques, more linear solvers, parallelizing the code, and adding to the \texttt{IRtools} \cite{GazzolaHansenNagy2018} package to include this algorithm.

The rest of the paper will proceed as follows. In section \ref{sec: BCD} we give a brief background to the CT problem, present the algorithm, and discuss some theory behind  the algorithm. Next in section \ref{sec: accel} we provide an alternate form of the algorithm as a fixed point iteration and discuss acceleration techniques to improve convergence. Section \ref{sec: implementation} outlines considerations for implementation including parallelization. Numerical experiments and results are outlined in section \ref{sec: results}. Future directions are presented in section \ref{sec: conclusion}. An example on how to set up and run an experiment is given in the appendix.

\section{Block Coordinate Descent}\label{sec: BCD} 

In this section, we discuss the mathematical concepts and modeling of the CT image reconstruction problem as well as different techniques involved to effectively solve the inverse, ill-posed problem. We begin by setting up the computed tomography problem which boils down to finding attenuation coefficients of an object made up of multiple materials. An image can then be constructed by a color mapping based on the attenuation coefficients. For an object made up of a single material, Beer's Law describes the amount of radiation that can pass through it \cite{Eps08}:
\begin{equation}\label{equa: Beer}
    I = I_0 e^{-\mu d}.
\end{equation}
In the above equation $I_0$ is the initial energy of the X-ray that goes into an object, and $I$ is the energy of the X-ray that leaves it. $d$ is the distance that the X-ray beam travels through the object and $\mu$ is the attenuation coefficient. Figure \ref{fig:Beer2} illustrates these parameters.

\begin{figure}[H]
    \centering
    \includegraphics[width=6cm]{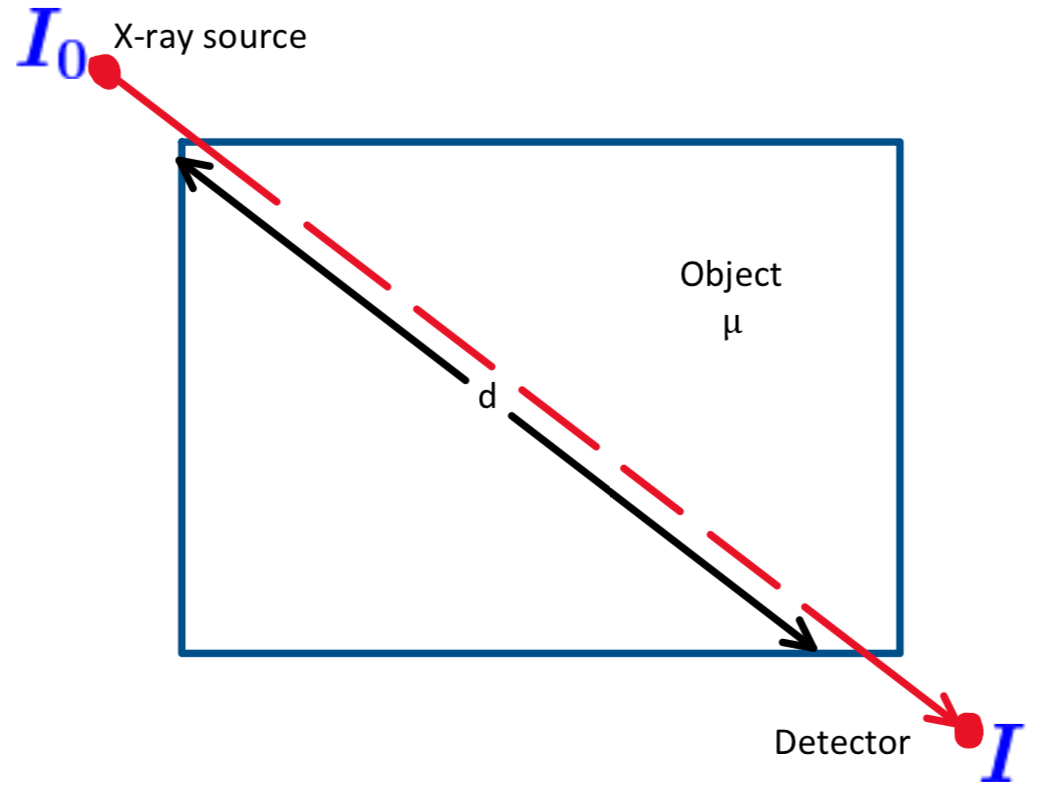}
    \caption{Beer's Law for a Single Material Object}
    \label{fig:Beer2}
\end{figure}
Beer's law implies that
$$
-\log\left(\dfrac{I}{I_0}\right) = \mu d.
$$
Thus, if the initial and final radiation amounts are known, as well as the distance $d$, finding the attenuation coefficient is as simple as solving a linear equation. In practice objects of interest are made up of several materials, which each have their own attenuation coefficients as illustrated in figure \ref{fig:Beer2}.

\begin{figure}[H]
    \centering
    \includegraphics[width=6cm]{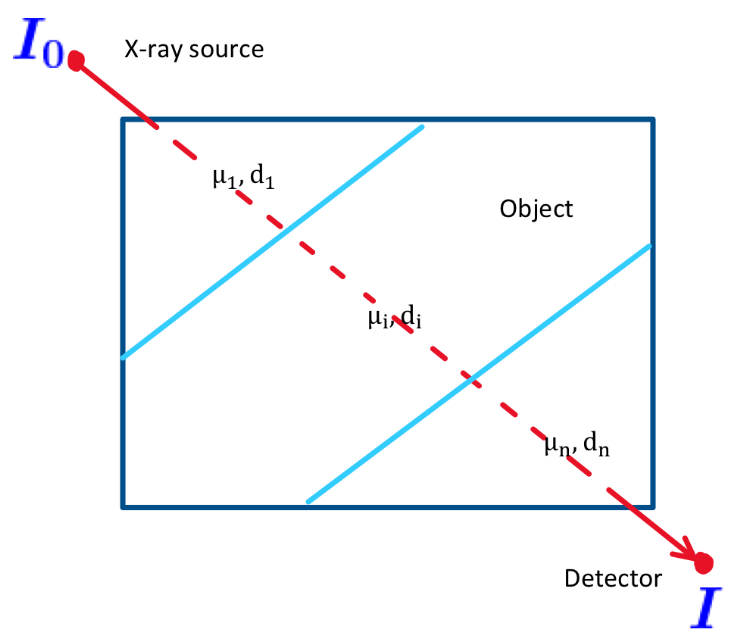}
    \caption{Beer's Law for an Object of Multiple Materials}
    \label{fig:Beer2}
\end{figure}
 Applying Beer's law to each different material, assuming there are $n$ materials, the X-ray energy leaving the object can be written as
 $$I = I_{0} e^{-\sum_{j=1}^n \mu_jd_{j}}$$
 which implies
 $$-\log\left(\dfrac{I}{I_{0}}\right) = \sum_{j=1}^n \mu_jd_{j}.$$
 This gives one linear equation in multiple variables, and so additional X-ray beams, at different angles, need to be transmitted through the object to obtain $n$ linearly independent equations for the $n$ variables.
 
\begin{figure}[H]
    \centering
    \includegraphics[width=6.5cm]{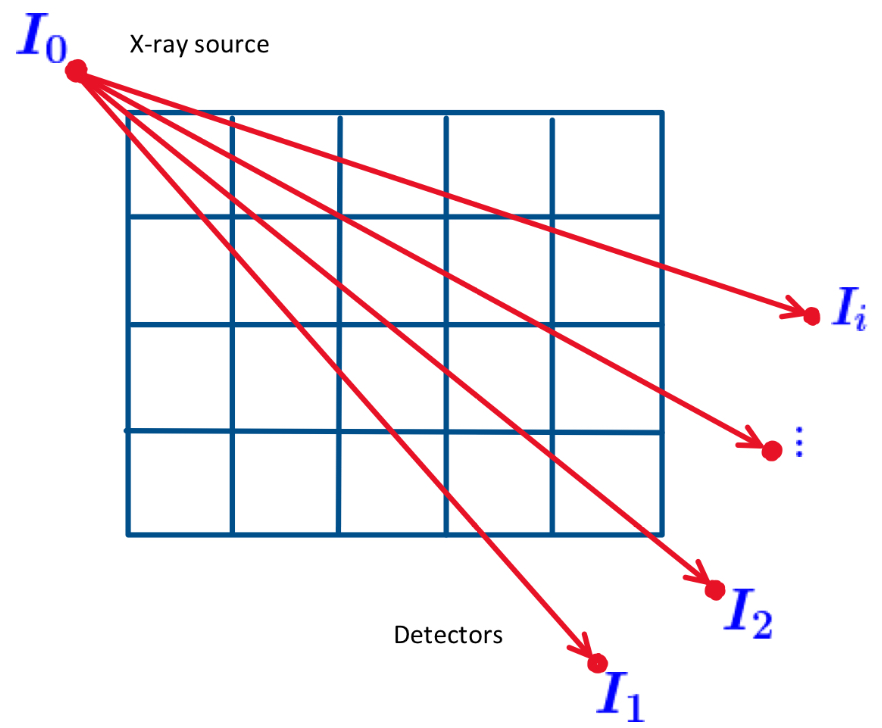}
    \caption{Fan-beam X-rays}
    \label{fig:Beer3}
\end{figure}

 Most often in practice it is not known precisely where each material begins and ends. To overcome this a pixel grid is overlaid on the object, and the attenuation of each pixel is sought. This involves taking measurements using many X-ray projections from different angles. Often in practice the X-ray source projects multiple X-rays spread out like a fan as illustrated in figure \ref{fig:Beer3}. In addition, the source and detector rotate around the object to obtain additional measurements.

Now suppose that there are $n$ pixels and $m$ X-ray beams are projected through the object. Let $d_{ij}$ be the distance that the $ith$ X-ray beam travels through pixel $j$. $I_0$ is still the initial energy, and let $I_i$ be the energy of the $i$th X-ray beam as it leaves the object and hits the detector. Then the relationship between the unknown attenuation coefficients and the known distances and energies can be written as a system of equations
\begin{equation}\label{equa: GenProb}
     \underbrace{\begin{bmatrix} d_{11} & d_{12} & \cdots & d_{1n} \\ d_{21} & d_{22} & \cdots & d_{2n} \\ \vdots & \vdots & \ddots & \vdots \\
    d_{m1} & d_{m2} & \cdots & d_{mn}\end{bmatrix}}_{\mb{A}}\underbrace{\begin{bmatrix} \mu_1 \\ \mu_2 \\ \vdots \\ \mu_n \end{bmatrix}}_{\mathbf{x}} = \underbrace{\begin{bmatrix}-\log\left({I_1 /I_{0}}\right)\\-\log({I_2 /I_{0}})\\\vdots\\-\log\left({I_m /I_{0}}\right)\end{bmatrix}}_{\mathbf{b}}.
\end{equation}

In this system of equations the vector $\mb{b}$ is known as the \textit{sinogram}. Note that the matrix $\mb{A}$ in the above equation will be sparse because each X-ray will only pass through a few of the pixels. 

Now, notice that each distance $d_{ij}$ is a function of the X-ray source location. Let $R$ denote the distance of the X-ray source to the center of the object, and $\theta$ denote the angle of rotation of the source relative to the center of the object. Then the matrix $\mb{A}$ is a function of the vector $\mb{p}$, denoted as $\mb{A}(\mb{p})$, where $\mb{p}$ contains the $R_i$, and $\theta_i$ for each source location (that is each time a fan of X-ray beams is emitted).


In a standard computed tomography problem with known geometry parameters, solving equation \eqref{equa: GenProb} requires regularization due to the ill-conditioning of $\mb{A}$ \cite{HansenJorgensenLionheart2021}. In this paper, the problem we seek to study has unknown geometry parameters, and can be stated as

\begin{equation}\label{eqn: CTprob}
 \mathbf{x} = \argmin_{x,p} ||\mb{A}(\mathbf{p})\mathbf{x} - \mathbf{b}||_2^2 + \alpha^2||\mb{x}||_2^2.   
\end{equation}

Here $\alpha$ is a regularization parameter, which allows us to get a more accurate solution when the matrix $\mb{A}$ is ill-conditioned. The regularization parameter recasts the ill-conditioned problem to a nearby well conditioned one. We refer the reader to \cite{Hansen},\cite{Bjorck} for information on the role of regularization in inverse and least squares problems. The solution to \eqref{eqn: CTprob} can be approximated by using an alternating minimization scheme known as block coordinate descent or BCD, which we describe in algorithm 1.

\begin{algorithm}[H]
\textbf{Input} $\mathbf{p_0} \in \R^{m}$\\
\textbf{Output} $\mathbf{x_k}$\\
  \For{$k = 0,1,2,...$}{
    $\mathbf{x_{k+1}} = \argmin ||\mb{A}(\mathbf{p_k})\mathbf{x} - \mathbf{b}||_2^2 + \alpha^2 ||\mathbf{x}||_2^2$\\
    $\mathbf{p_{k+1}} = \argmin ||\mb{A}(\mathbf{p})\mathbf{x_{k+1}} - \mathbf{b}||_2^2$\\
    }
 \caption{Tomography Block Coordinate Descent}
 \label{alg: BCD}
\end{algorithm}

Thus, equation \eqref{eqn: CTprob} can be approximated by an alternating minimization scheme where one iteration involves solving a large linear least squares problem, and a large non-linear least squares problem. To solve the large regularized linear least squares problem in step 4, we mainly use a hybrid conjugate gradient method (LSQR) that can automatically choose regularization parameters. Mathematical details describing this method can be found in \cite{WGCV}. The implementation we use is described in \cite{GazzolaHansenNagy2018} (specifically, we use the method called IRhybrid\_lsqr). A comparision of this to other linear least squares solvers is given in section \ref{sec: lls}. To solve the non-linear least squares problem in step 5 we use a projected quasi-Newton method in Kelley's implicit filtering package \cite{Kelley2011}.


\section{Acceleration Techniques}\label{sec: accel}
In this section we discuss techniques to accelerate the convergence of the BCD algorithm. The key is to recast the problem in terms of a fixed point problem.

\begin{algorithm}[H]
\textbf{Input} $\mb{x}_0 \in \R^{n}$\\
\textbf{Output} $\mb{x}_k$\\
  \For{$k =0, 1,2,...$}{
    $\mb{p}_{k+1} = \argmin ||\mb{A}(\mb{p})\mb{x}_{k} - \mb{b}||_2^2$\\
    $\mb{x}_{k+1} = \argmin ||\mb{A}(\mb{p}_{k+1})\mb{x} - \mb{b}||_2^2 + \alpha^2 ||\mb{x}||_2^2$\\
    }
 \caption{Fixed Point Tomography Reconstruction}
 \label{alg: fixed point}
\end{algorithm}

 Now we define a function $g$ that is one iteration of algorithm \ref{alg: fixed point}. That is assuming the $\mathbf{b}$, $\alpha$, and the function $\mb{A}(\mb{p})$ are known define $g(\mb{x}_k)$ as
 $$g(\mb{x}_k) := \argmin ||\mb{A}(\mb{p}_{k+1})\mb{x} - \mb{b}||_2^2 + \alpha^2 ||\mb{x}||_2^2$$
 where $\mb{p}_{k+1}$ is defined as
 $$\mb{p}_{k+1} := \argmin ||\mb{A}(\mb{p})\mb{x}_{k} - \mb{b}||_2^2.$$
 
 Thus, the problem of finding an approximate solution of \eqref{eqn: CTprob} becomes that of finding a fixed point of $g$. This is useful as acceleration techniques for fixed point problems have been widely studied \cite{RAMIERE2015}. It should be noted that the problem can also be stated in terms of finding a fixed point of $g(\omega_k)$ where 
  $$g(\omega_k) := \argmin ||\mb{A}(\mb{p})\mb{x} - \mb{b}||_2^2 + \alpha^2 ||\mb{x}||_2^2$$
 and $\omega_k := (\mb{x}_k,\mb{p}_k)$. However numerically this performed worse as shown in subsection \ref{subsec: accelerationTest}.

 In our numerical experiments we used three fixed point acceleration schemes, and the remainder of the section will be devoted to discussing them. For the remainder of the section we will denote a function as $f(x)$ if it is a function is a single variable, and $F(\mb{x})$ if it is a function of several variables.

We motivate the first two acceleration methods starting with the single variable fixed point problem of finding a point $x_k$ such that $f(x_k) = x_{k+1}$. To begin, we let $\Delta$ denote the difference operator, that is let $\Delta x_k := f(x_{k}) - x_k$, $\Delta f(x_k) := f(f(x_k)) - f(x_k)$, and $\Delta^2 x_k := \Delta f(x_k) - \Delta x_k$. Aitken's $\Delta^2$ method \cite{Aitkens} along with its recursive application and the second order Steffensen method are popular choices for solving single variable fixed point problems \cite{Steffensen}. The Aitken $\Delta^2$ method is derived from approximating the multiplying constant for a linearly converging sequence, and can be stated as
$$x_{k+1} = x_k - \frac{(f(x_n) - x_n)^2}{f(f(x_k)) - 2f(x_k) + x_k} = x_k - \frac{(\Delta x_k)^2}{\Delta^2x_k}.$$
Or equivalently 
\begin{equation}\label{eqn: deltasquared}
    x_{k+1} = f(f(x_k)) - \frac{(\Delta f(x_k))^2}{\Delta^2x_k}.
\end{equation}
Vector generalizations of this algorithm have been extensively studied, see \cite{RAMIERE2015} for several references. Most come from defining the inverse of a vector $\mathbf{x}$ as 
$$\mb{x}\inv := \mb{x}\frac{1}{||\mb{x}||_2^2}.$$
Application of this definition to \eqref{eqn: deltasquared} leads to the \textit{Irons-Tuck} method \cite{IT}
$$\mb{x}_{k+1} = F(F(\mb{x}_k)) - \frac{\Delta F(\mb{x}_k) \cdot \Delta^2 \mb{x}_k}{|| \Delta^2 \mb{x}_k||^2}\Delta F(\mb{x}_k).$$
Where $\cdot$ denotes the dot product, $F(\mb{x}_k)$ is the vector valued function we seek to find a fixed point of, and the $\Delta$ operator is defined as in the scalar case. Numerically the Irons-Tuck method has been shown to outperform many other vector generalizations of the $\Delta^2$ process \cite{Compare}, but is computationally expensive when $F(\mb{x}_k)$ is expensive to evaluate. One alternative is to use the \textit{crossed secant} method
$$\mb{x}_{k+1} = F(\mb{x}_k) - \frac{(F(\mb{x}_k) - F(\mb{x}_{k-1})\cdot (\Delta \mb{x}_k - \Delta \mb{x}_{k - 1})}{||(\Delta \mb{x}_k - \Delta \mb{x}_{k - 1})||_2^2} \Delta \mb{x}_k.$$
It has the advantage of only needing one evaluation of $F$ per iteration, and in many numerical tests performs similar to the Irons-Tuck method \cite{RAMIERE2015}. In fact, if each iteration of the crossed secant method is alternated with a standard fixed point iteration the resulting method is the Irons-Tuck method.

Another alternative is \textit{Anderson Acceleration}, sometimes called \textit{Anderson Mixing}. The general idea is to use information from a predetermined number of residuals to increase convergence. Here we only state a more easily implementable version of the algorithm. For a more general algorithm see \cite{AA}. Let $r$ be the number of residuals to use, and denote
$$\mathcal{X}_k = (\Delta \mb{x}_{k-r} , \cdots, \Delta \mb{x}_{k - 1}).$$
Then the Anderson Acceleration has the following form:

\begin{algorithm}[H]
\textbf{Input} $\mb{x}_0 \in \R^{n}$, $r \ge 1$\\
\textbf{Output} $\mb{x}_k$\\
Compute $F(\mb{x}_0) = \mb{x}_0$
  \For{$k = 1,2,...$}{
    $r_k = \min\{r,k\}$\\
    $\gamma^{(k)} = \argmin ||\Delta \mb{x}_{k} - \mathcal{X}_k \gamma||_2$\\
    $\mb{x}_{k+1} = F(\mb{x}_k) - \sum_{i = 1}^{r_k}\gamma^{(k)}_i\left[(F(\mb{x}_{k - r_k +i + 1}) - F(\mb{x}_{k - r_k + i})\right]$
}
 \caption{Anderson Acceleration}
\end{algorithm}
Implementation details can be found in \cite{AAImplementation}. In short, the linear least squares problem can be efficiently solved by updating a $QR$ factorization of $\mathcal{X}$ after every iteration. Additionally, we implemented a hyperparameter which removes a column of $\mathcal{X}$ if the condition number of the upper triangular matrix in the $QR$ factorization becomes too large. In section \ref{subsec: accelerationTest} we compare the use of these three algorithms in our implementation.

\section{Implementation}\label{sec: implementation}
In this section we discuss practical considerations for the implementation of algorithm \ref{alg: BCD}. Recall that each time the X-ray source fires it is associated to one angle and one $R$ parameter. Therefore $\mb{A}$ is created as 
$$\mb{A}(\mb{p}) = \begin{bmatrix}\mb{A}(\theta_1,R_1)\\ \vdots \\ \mb{A}(\theta_n,R_n) \end{bmatrix}.$$

 The non-linear least squares problem can then be solved by solving subproblems of the form
$$ \underset{\theta_i,R_i}{\min} ||\mb{A}(\theta_i,R_i)\mb{x}_{k} - \mb{b}_i||_2^2$$
where each $\mb{b}_i$ is the corresponding part of the $\mathbf{b}$ vector. This is advantageous as it allows for easy parallelization of algorithm \ref{alg: BCD} leading to significant speedup, and a comparison is made in subsection \ref{subsec: parallel test}. 

To simulate the problem we will use the IRTools package \cite{GazzolaHansenNagy2018}. IRTools provides several state of the art regularized linear least squares solvers which we will use in our implementation of algorithm \ref{alg: BCD}. We will compare three linear least squares solvers available in the IRtools package which use different forms of regularization. The first is the hybrid-LSQR algorithm \cite{WGCV}, which regularizes by adding a two-norm penalty term. The next is the iteratively re-weighted norm approach or IRN \cite{rodriguez-2008-sparse}, which uses a one-norm regularization. Finally, the fast iterative shrinkage threshold algorithm or FISTA \cite{fista} is used which has no regularization term, but constrains the solution vector $\mb{x}$ such that $\mb{x} \ge 0$.  A comparison of these methods is given in subsection \ref{sec: lls}.

To approximate the solution of the non-linear least squares problem  we compared two methods. The first is a  Trust-Region-Reflective algorithm studied in \cite{lsq2},\cite{lsq1} and implemented in the MATLAB function \texttt{lsqnonlin}. The second is the method of implicit filtering \cite{Kelley2011} which uses a projected quasi-Newton iteration using difference gradients. Both methods are for bounded problems with unknown derivative information, making them suitable for our problem, as the geometry parameters can be realistically bounded. In our numerical tests we found that implicit filtering led to a smaller image error and therefore it was used for all numerical experiments in the next section.

\section{Numerical Results}\label{sec: results}
In this section we highlight several experiments for the tomography problem. In each problem an initial guess for the $R$ parameters and angle parameters were given. Then a perturbation generated uniformly was added to each angle and $R$ value, and the BCD was run to produce a better image. For all tests Gaussian noise was added to the sinogram $\mb{b}$ creating a new vector $\mathbf{b}_{noise} = \mathbf{b} + \mathbf{\eta}$, where $||\mathbf{\eta}||_2 / ||\mathbf{b}||_2$ is equal to a specified noise level, in this case $0.01$.
The initial guess for the view angles was always $0:2:358$, and the initial guess for the $R$ parameters was $2$.  We will use the abbreviation $U(a,b)$ for a uniform distribution when describing perturbations to $R$ and angle parameters. In the graphs below the angle and $R$ parameter error represent the relative error between the true perturbations and the current approximation of those perturbations. The $\ell_2$-norm was used for relative errors in each case and an initial guess of $0$ was given for all perturbations. 

The images used were the Shepp-Logan phantom (Figure \ref{fig:iterbraincomp}a), and the spine image from the MATLAB image-processing toolbox (Figure \ref{fig:AccelSpines}: True Image), with each image size being $256 \times 256$ unless otherwise stated. All tests were run on a Microway system that has four Intel Xeon E5-4627 CPUs with 40 cores and 1 TB of memory running Ubuntu Linux. An example of setting up and running a test is given in appendix \ref{sec: Code Example}. In the rest of this section we will first show a scaling study demonstrating the effectiveness of parallelizing the code. We will then show that the fixed-point acceleration techniques produce a smaller image error and converge faster than not using an acceleration technique. Finally, we compare several solvers for the linear least squares solver.

\subsection{Parallel Speedup}\label{subsec: parallel test}

In this section we demonstrate the effectiveness of the parallelization discussed in section \ref{sec: implementation}. Table \ref{tab: parallel times} shows a 20 iteration run of the same simulation using different image sizes with the image size being $N \times N$ where $n = N^2$. In all the timing tests the Shepp-Logan phantom was used with perturbations as realizations of the distribution $U(-0.25,0.25)$ added to each angle and $R$ value. The hybrid-LSQR algorithm was used as the linear least squares solver. The aforementioned computer was used with twelve workers.

\begin{table}[H]
    \centering
    \begin{tabular}{|c|c|c|c|c|c|}
    \hline
         N        & 32 & 64 & 128 & 256 & 512\\
    \hline
         Serial   & 257 & 714 & 1308 & 4238 &  13771 \\
    \hline
         Parallel & 35 & 166 & 289 & 1026 & 3810\\
    \hline
    \end{tabular}
    \caption{Run Times in Seconds} 
    \label{tab: parallel times}
\end{table}

From the results in table \ref{tab: parallel times} we see drastic speedup in all cases with the parallelization. For larger problems solving the linear least squares problem appears to take most of the time, since $\mb{A}$ is a large matrix. For example when $N$ is 256 the size of $\mb{A}$ is $65160 \times 65536$ and when $N$ is 512 the size of $\mb{A}$ is $131400 \times 266256$. Regardless, the parallel implementation made possible by the problem model provides significant speedup. 

\subsection{Acceleration Test}\label{subsec: accelerationTest}

\begin{figure}[H]
\begin{center}
\begin{tabular}{ccc}
\includegraphics[width=5.5cm]{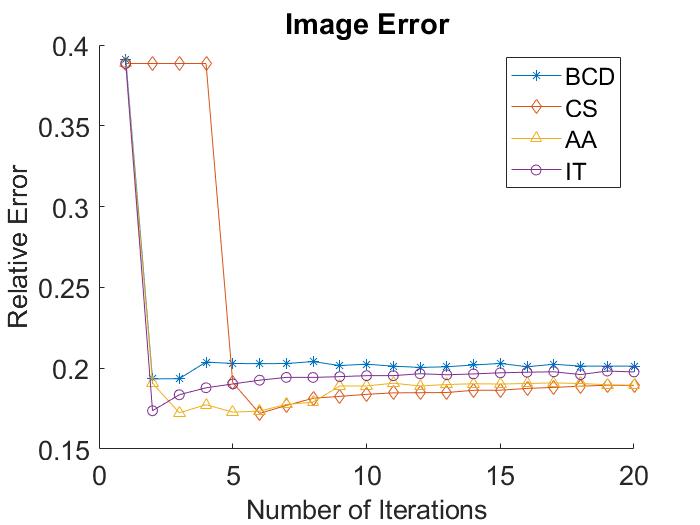} &
\includegraphics[width=5.5cm]{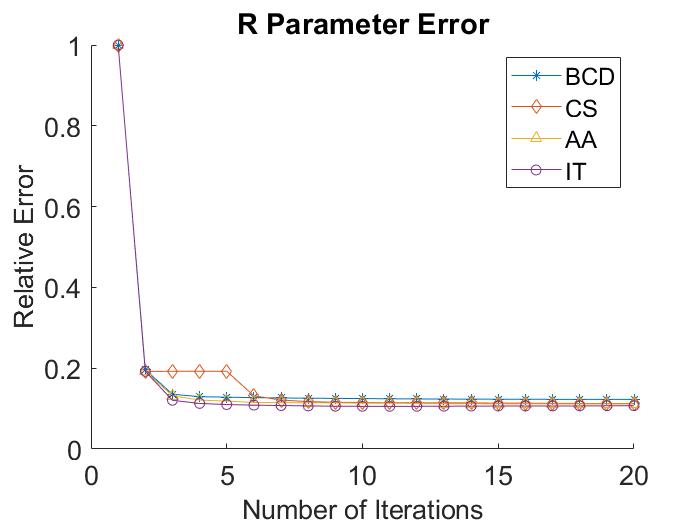} 
\includegraphics[width=5.5cm]{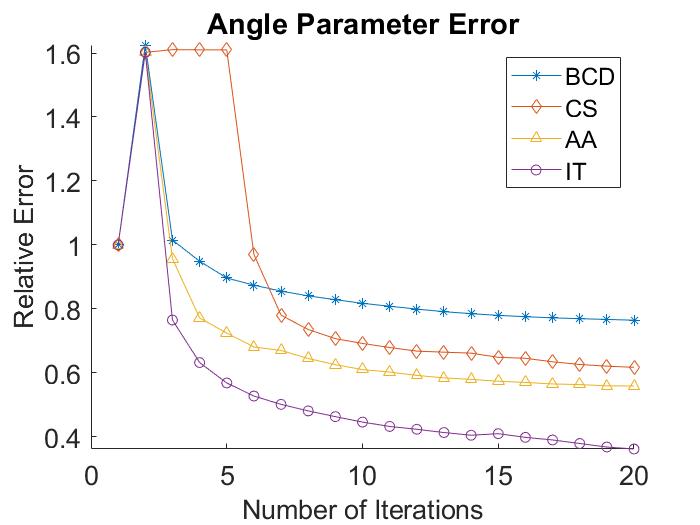} 
\end{tabular}
\end{center}
\caption{Graphs Comparing Acceleration Techniques} \label{fig:AccelComp}
\end{figure}

The next test we ran was to compare the different acceleration techniques using the spine image. Perturbations were added as realizations of $U(-0.5,0.5)$. The hybrid-LSQR algorithm was used as the linear least squares solver, and the accelerated BCD algorithms were run for twenty iterations. In the legends in figure \ref{fig:AccelComp}, BCD is for no acceleration, CS for crossed secant, AA for Anderson Acceleration, and IT for the Irons-Tuck method. To give reference to table \ref{tab:Accel RunTimes} the tests were run on 12 local workers.

\begin{table}[H]
    \centering
    \begin{tabular}{|c|c|c|c|}
        \hline
          BCD & CS &AA & IT    \\
         \hline
         1004 & 1042 & 1095 & 2375\\
        \hline 
    \end{tabular}
    \caption{Run Times in Seconds}
    \label{tab:Accel RunTimes}
\end{table}

\begin{figure}[H]
    \includegraphics[width=.33\textwidth]{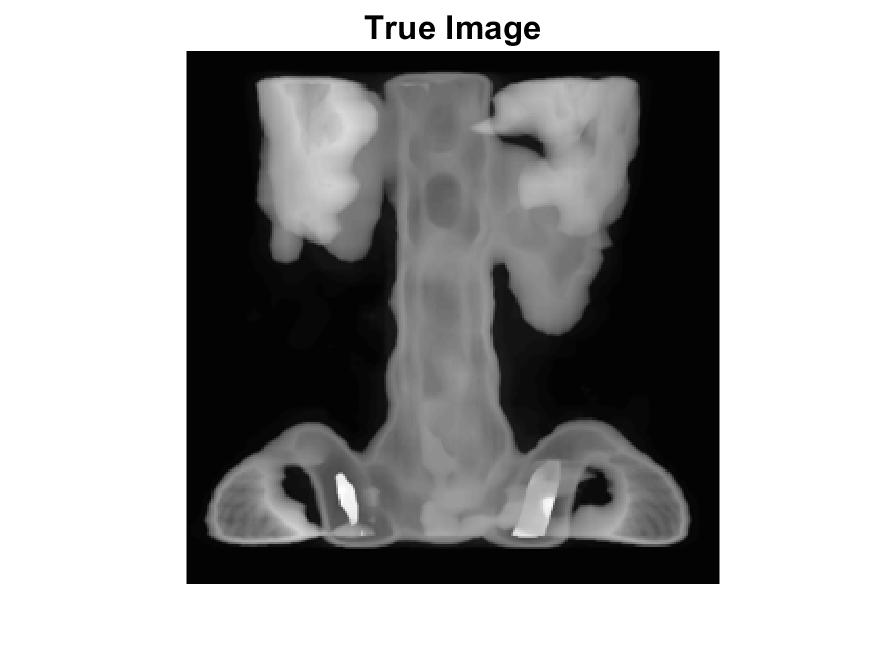}\hfill
    \includegraphics[width=.33\textwidth]{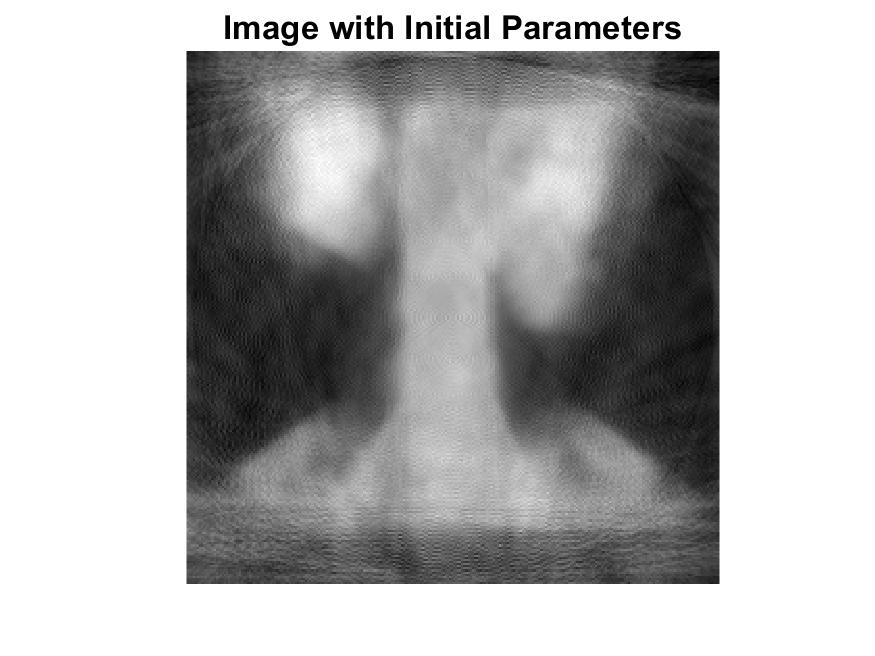}\hfill
    \includegraphics[width=.33\textwidth]{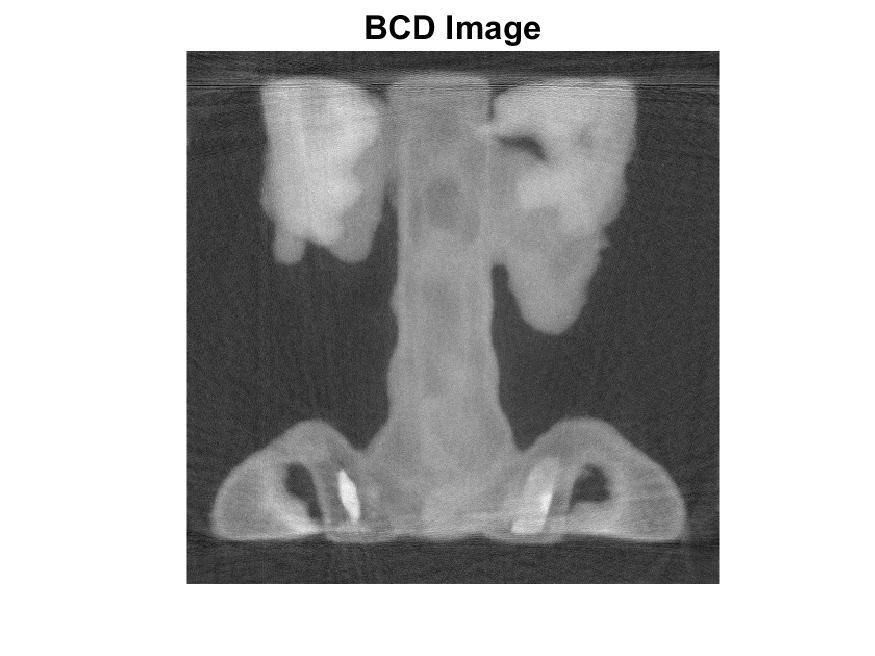}
    \\[\smallskipamount]
    \includegraphics[width=.33\textwidth]{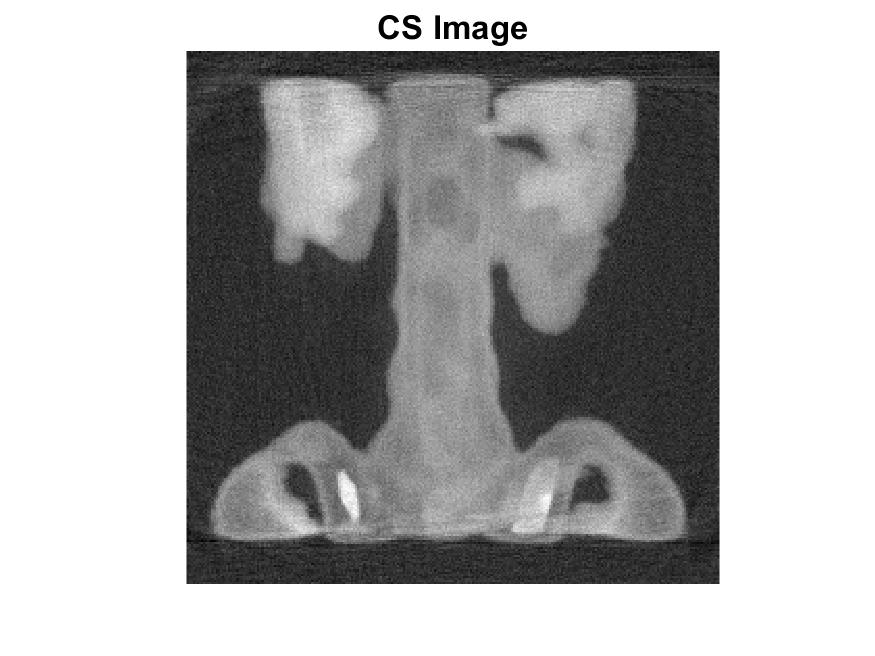}\hfill
    \includegraphics[width=.33\textwidth]{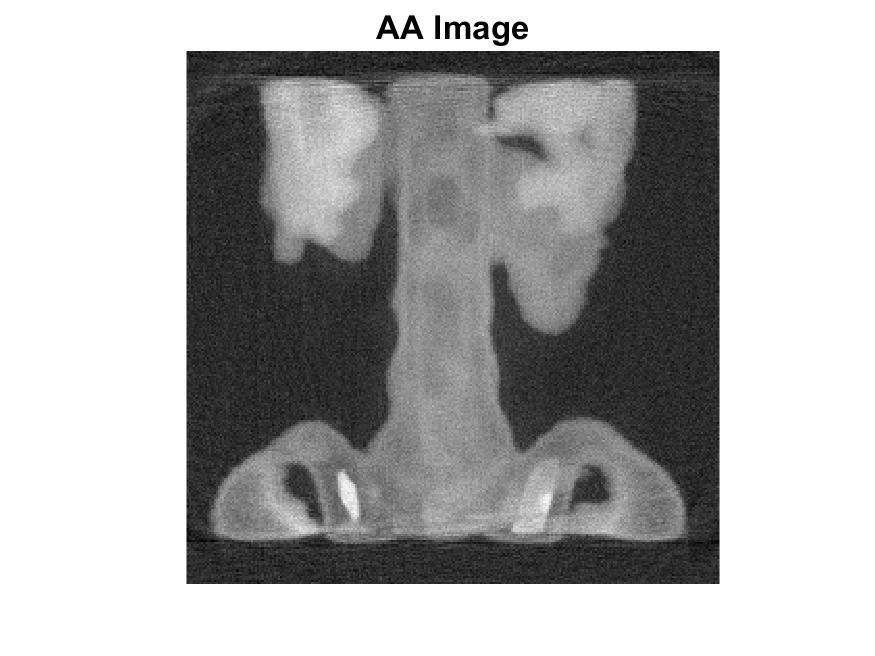}\hfill
    \includegraphics[width=.33\textwidth]{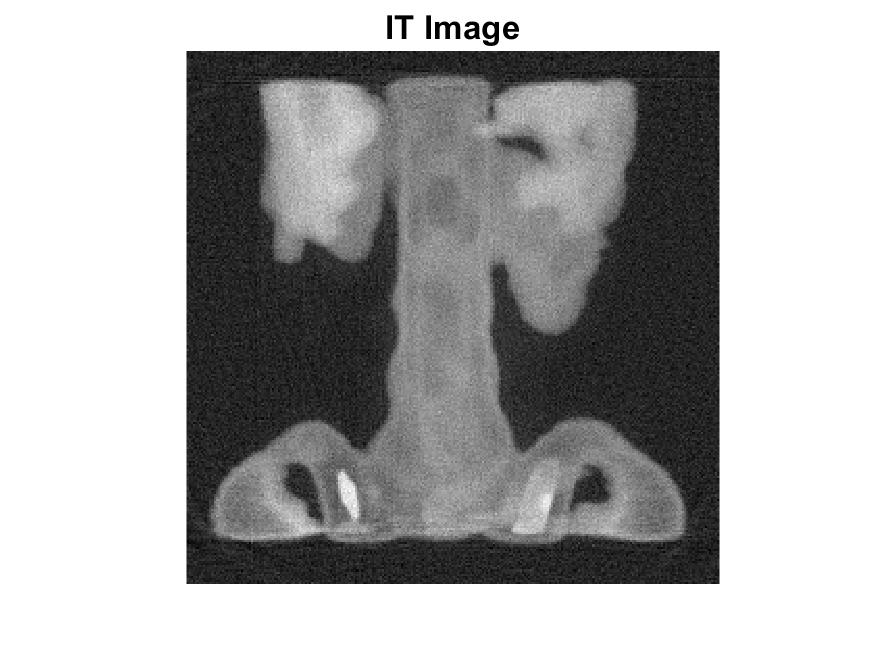}

    \caption{Spine Images after Parameter Optimization}
    \label{fig:AccelSpines}
\end{figure}

From these tests we see that the acceleration techniques converge to slightly smaller error norms, with the crossed secant method and Anderson Acceleration performing the best. The Irons-Tuck method converged much better in the angle parameters, but took much longer. This was expected as the function evaluation is quite expensive, and as previously noted Irons-Tuck requires an extra function evaluation at each iteration. Figure \ref{fig:AccelSpines} shows the true image, and the image after the parameter optimization. Despite having slightly larger image error, the Irons-Tuck image seems to have the least background noise.

As previously mentioned, in these tests we solved for the fixed point where the fixed point is the image vector. The table below compares the case when the fixed point vector is defined to be $\omega := (\mathbf{x}, \mathbf{p})$, with the same parameters as before.
\begin{table}[H]
    \centering
    \begin{tabular}{|c|c|c|c|c|}
        \hline
          Acceleration Scheme & Fixed Point Type & Angle Error & R Error &  Image Error  \\
         \hline
         Anderson & & & &  \\ 
         
            &    $\mathbf{x}$ & $0.557$ &$0.111$ & $0.189$ \\
            \hline 
            &   $\omega$   & $14.642$ & $14.020$ & $0.190$\\
        \hline
         Crossed Secant & & & &\\
                    & $\mathbf{x}$  & $0.614$ &$0.112$ &$0.189$\\
                    \hline
                    & $\omega$ & $15.251 $ &$14.079 $ & $0.1774$\\
                    \hline
         Irons-Tuck & & & & \\
         & $\mathbf{x}$ & $0.359$ & $0.108$ & $0.199$\\
         \hline
        & $\omega$ &   $14.009$   & $14.009$     & $0.195$\\
        \hline
    \end{tabular}
    \caption{Comparison of Fixed Point Problems}
    \label{tab:Fix the Point}
\end{table}

As can be seen from table \ref{tab:Fix the Point} when the $\omega$ fixed point scheme is used the image error remains approximately the same, while the error in the angle and $R$ parameters grows large. Thus, in our scheme we choose to only look for a fixed point of $\mb{x}$. The conclusion of this experiment overall is that when an accelerated fixed point scheme on $\mb{x}$ is used it may produce a smaller image error norm.

Having shown acceleration techniques looking for a fixed point of $\mb{x}$ to be effective in reducing error, we now demonstrate that giving a stopping tolerance they can effectively reduce the number of iterations needed to find a fixed point. Different from the other experiments the perturbations used here are sampled from  $U(-0.25,0.25)$ with other factors such as $R$ and view angle initial guesses, and image size held the same. We use a relative error stopping criterion of 
\begin{equation*}\label{eqn:Stop}
    \dfrac{||\mathbf{x}_k - \mathbf{x}_{k-1}||_2}{||\mathbf{x}_{k-1}||_2} \leq 0.03.
\end{equation*}

Below we show the results using BCD without acceleration and with Anderson Acceleration.

\begin{figure}[H]
    \begin{center}
        \begin{tabular}{cc}
            \includegraphics[width = 6.5cm]{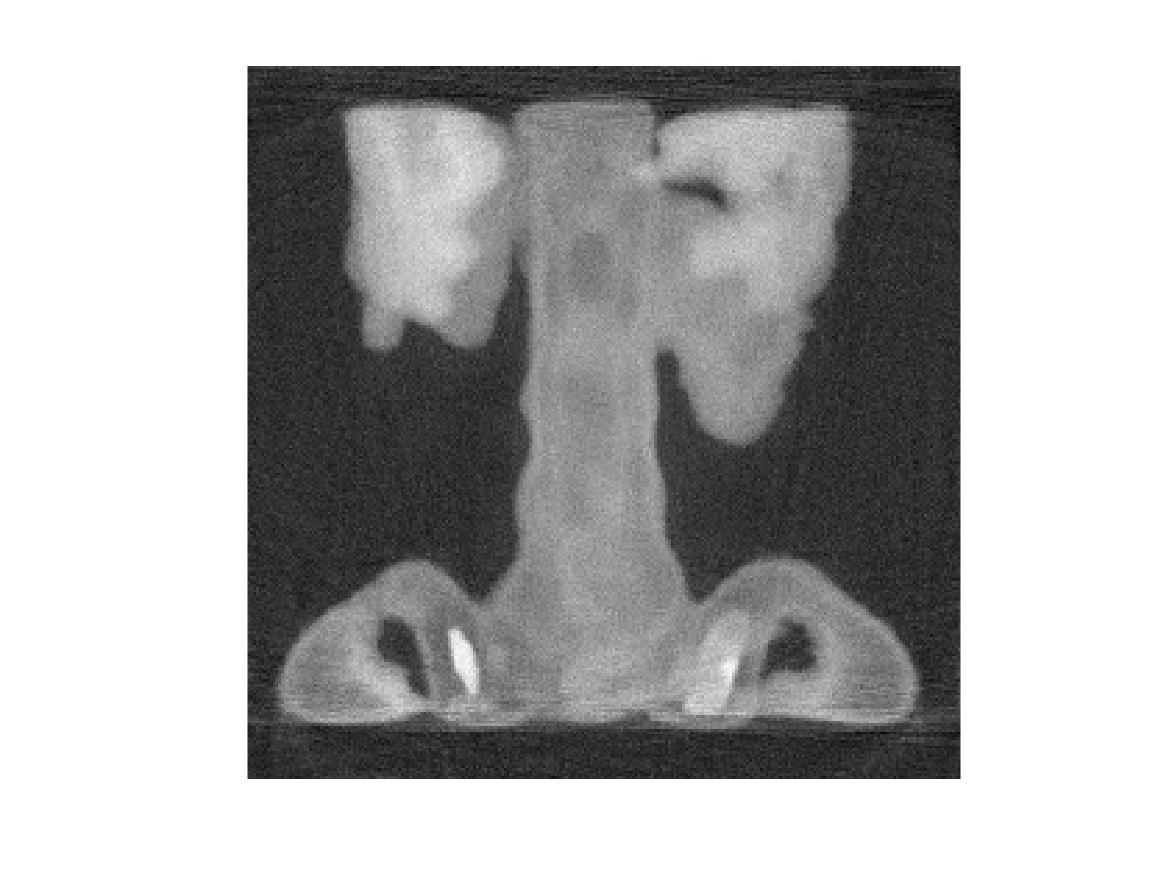}&
            \includegraphics[width = 6.5cm]{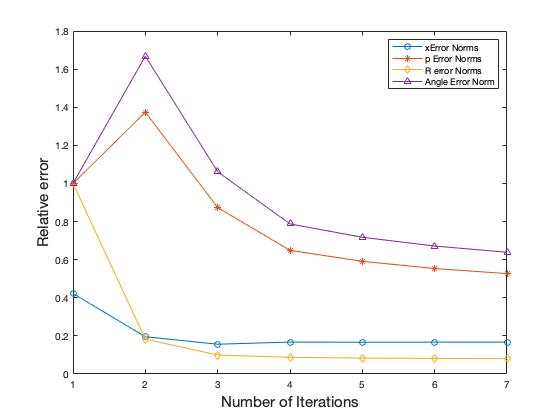}\\
            {\footnotesize Final Image}&
            {\footnotesize Error Plot}
        \end{tabular}
    \end{center}
\caption{Convergence of BCD without acceleration} \label{fig:BCDStop}
\end{figure}

\begin{figure}[H]
    \begin{center}
        \begin{tabular}{cc}
            \includegraphics[width = 6.5cm]{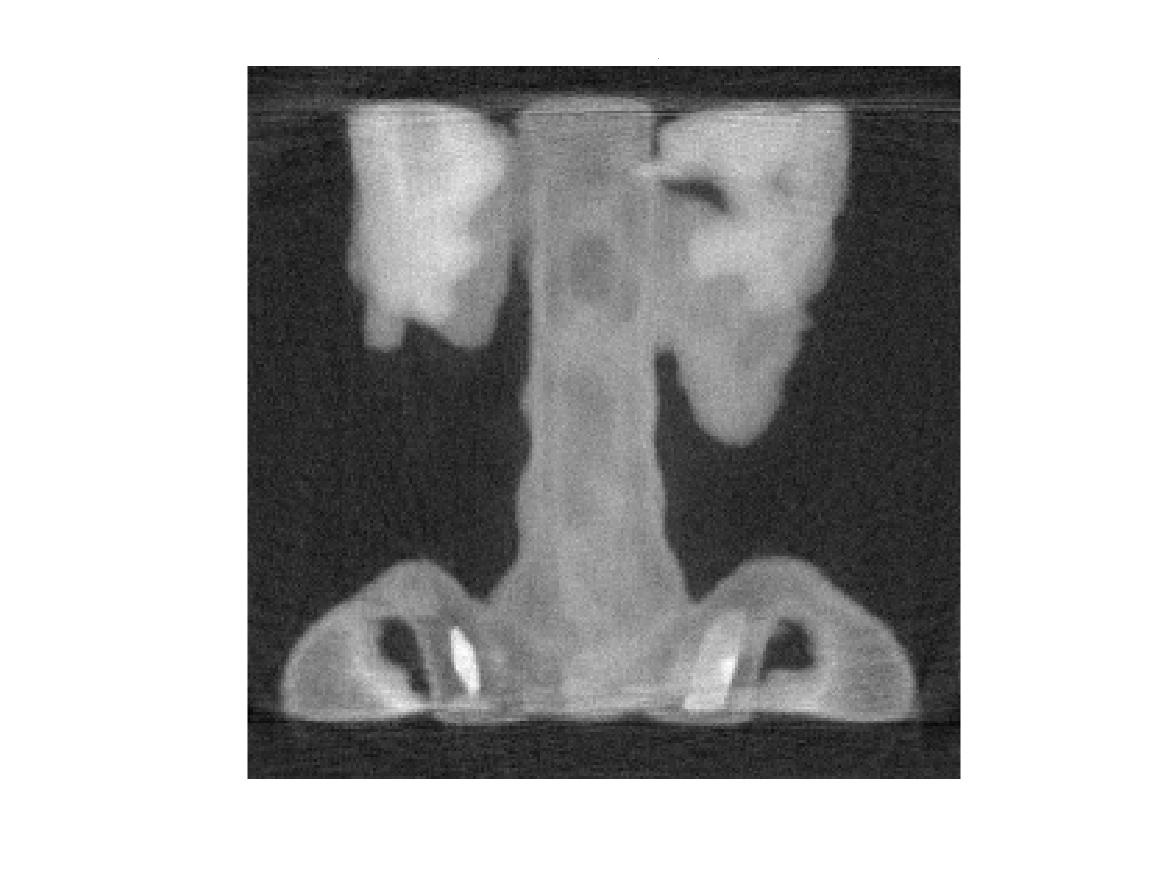}&
            \includegraphics[width = 6.5cm]{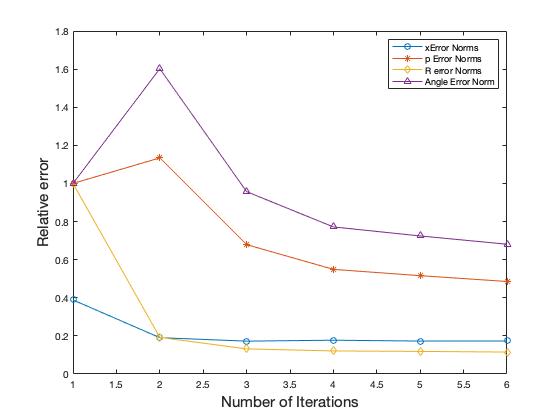}\\
            {\footnotesize Final Image}&
            {\footnotesize Error Plot}
        \end{tabular}
    \end{center}
\caption{Convergence of BCD with Anderson Acceleration} \label{fig:AAStop}
\end{figure}

 For this small example, using Anderson Acceleration caused the BCD to reach the stopping tolerance one iteration faster than without it. Since the fixed point function $g$ is expensive to evaluate saving one iteration reduces computation time non-trivially. For example, in the test performed in figure \ref{fig:AccelComp} one evaluation of $g$ averaged $50.2$ seconds to evaluate. 

\subsection{Comparison of Linear Least Squares Solvers}\label{sec: lls}
In this section we compare the performance of three linear least squares solvers available in the IRtools package using hybrid LSQR, FISTA, IRN. We did two tests using the Shepp-Logan phantom. Perturbations were chosen from realizations of $ U(-0.5,0.5)$. No acceleration techniques were used in these tests.

\begin{figure}[H]
\begin{center}
\begin{tabular}{ccc}
\includegraphics[width=5.5cm]{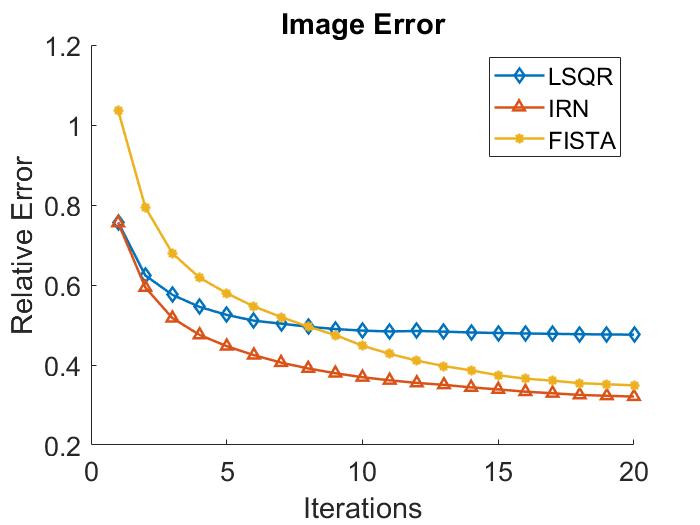} &
\includegraphics[width=5.5cm]{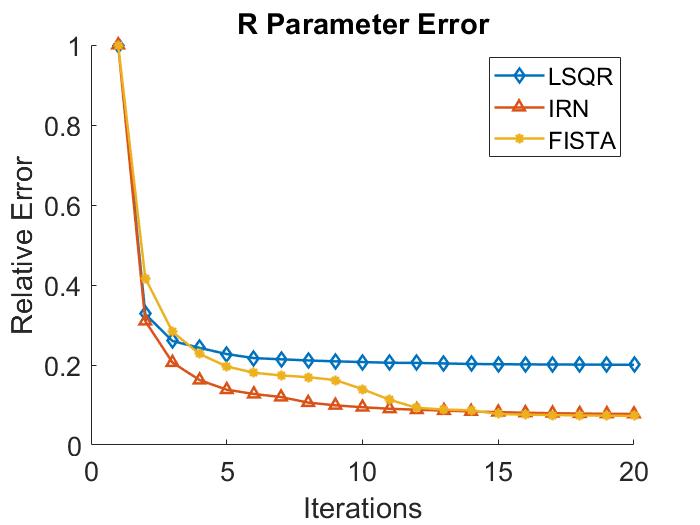} 
\includegraphics[width=5.5cm]{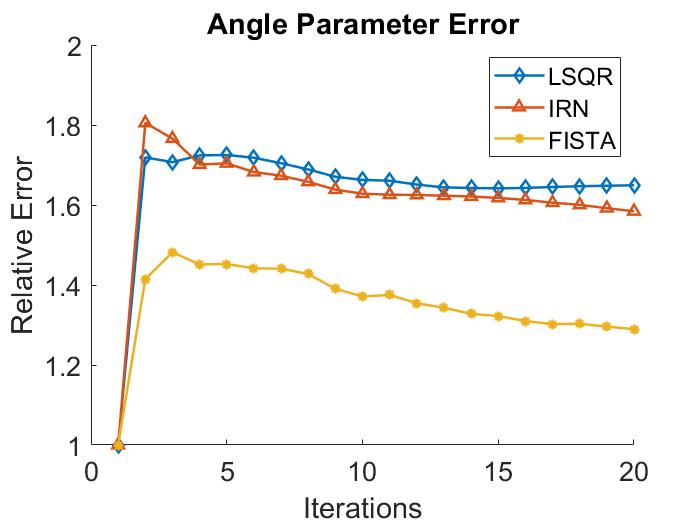} 
\end{tabular}
\end{center}
\caption{Graphs Comparing Linear Least Squares Solvers} \label{fig:IterComp1}
\end{figure} 

\begin{figure}[H]
\begin{center}
\begin{tabular}{cccc}
\includegraphics[width = 4cm]{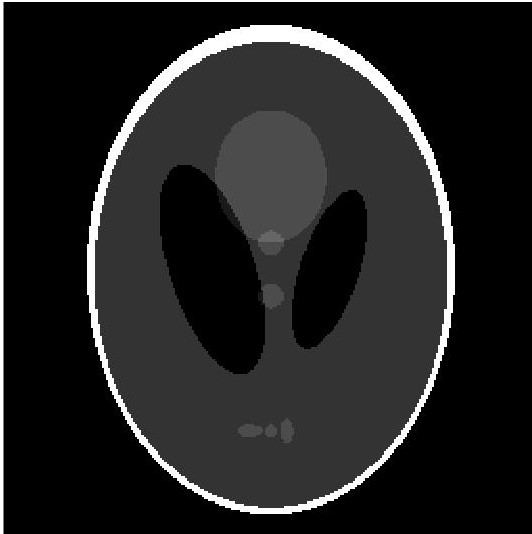}&
\includegraphics[width = 4cm]{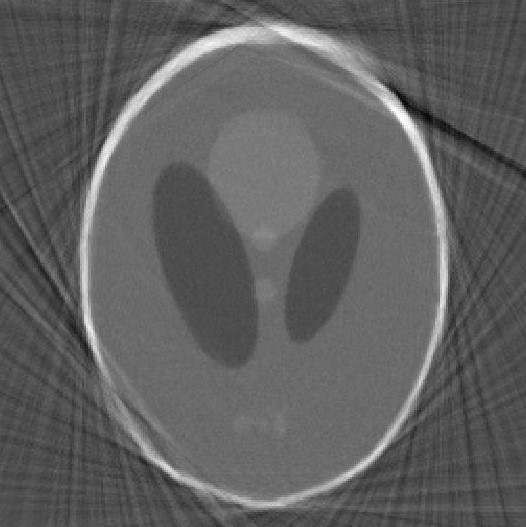}&
\includegraphics[width = 4cm]{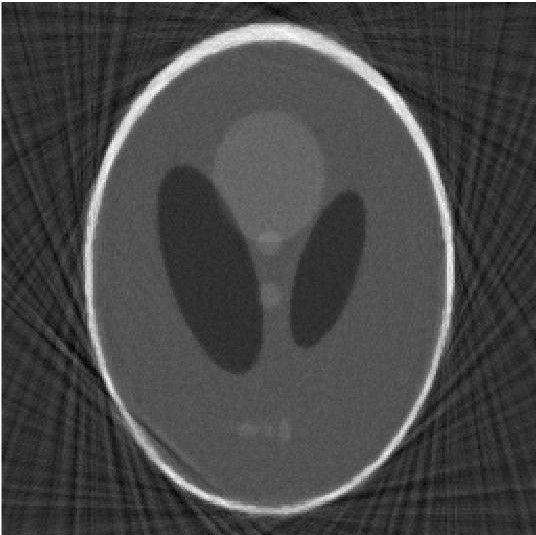}&
\includegraphics[width = 4cm]{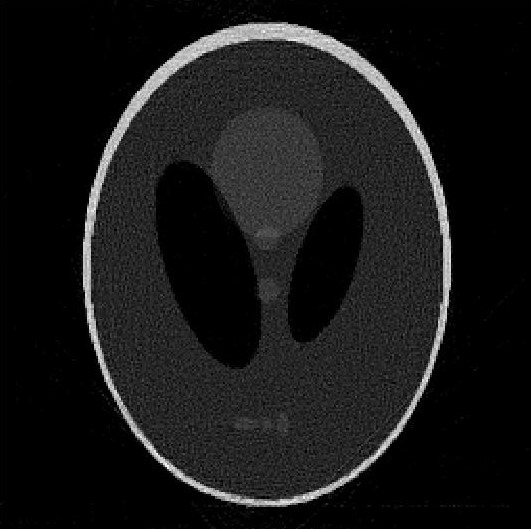}\\
{\footnotesize (a) True Solution} &
{\footnotesize (b) LSQR} &
{\footnotesize (c) IRN}&
{\footnotesize (d) FISTA}
\end{tabular}
\end{center}
\caption{Images for Linear Least Squares Comparison} \label{fig:iterbraincomp}
\end{figure}

In figure \ref{fig:IterComp1} we see that IRN and FISTA perform similarly, and drastically outperform the hyrbid-LSQR method. Interestingly all three methods perform poorly on the angle parameter estimation. Despite having similar image error norms, figure \ref{fig:iterbraincomp} illustrates that visually FISTA appears better than IRN.

\begin{figure}[H]
\begin{center}
\begin{tabular}{ccc}
\includegraphics[width=5.5cm]{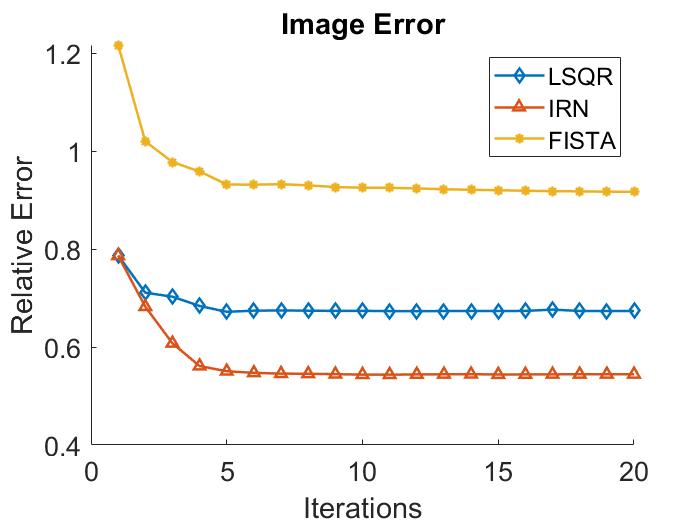} &
\includegraphics[width=5.5cm]{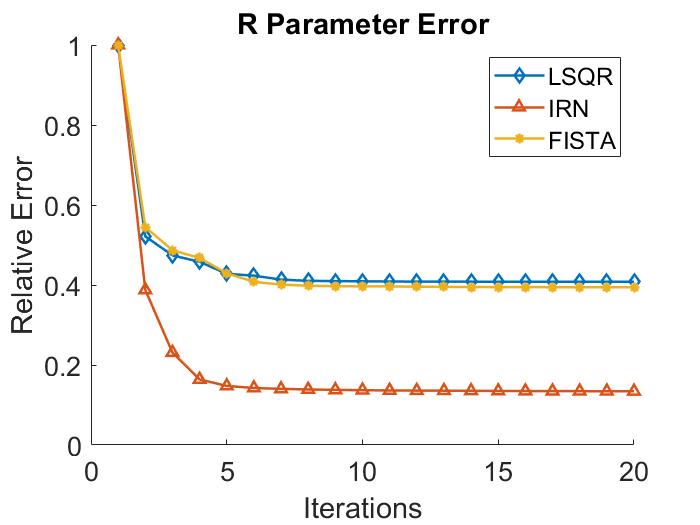} 
\includegraphics[width=5.5cm]{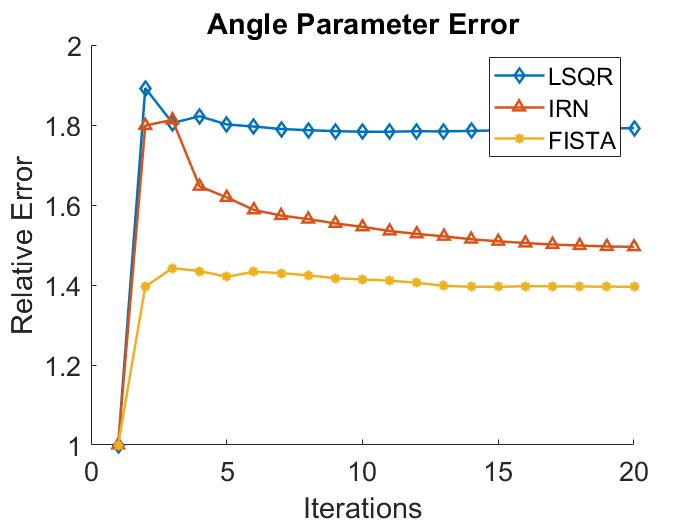} 
\end{tabular}
\end{center}
\caption{Graphs Comparing Linear Least Squares Solvers For Second Test Problem} \label{fig:IterComp2}
\end{figure}

The second test comparing the linear least squares solvers uses the same parameters as the previous test, except the perturbations are samples of $U(-1,1)$. This is a significant amount of perturbation in this case as it is up to half the original value of $R$ and half the distance between each angle. Figure \ref{fig:IterComp2} shows the relative errors.

\begin{figure}[H]
\begin{center}
\begin{tabular}{cccc}
\includegraphics[width = 4cm]{Images/IterTrue1.jpg}&
\includegraphics[width = 4cm]{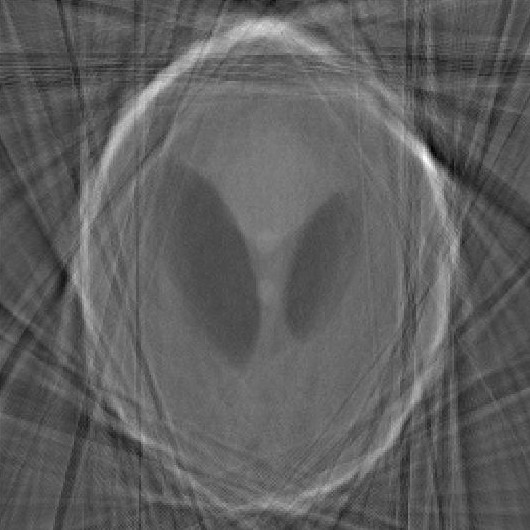}&
\includegraphics[width = 4cm]{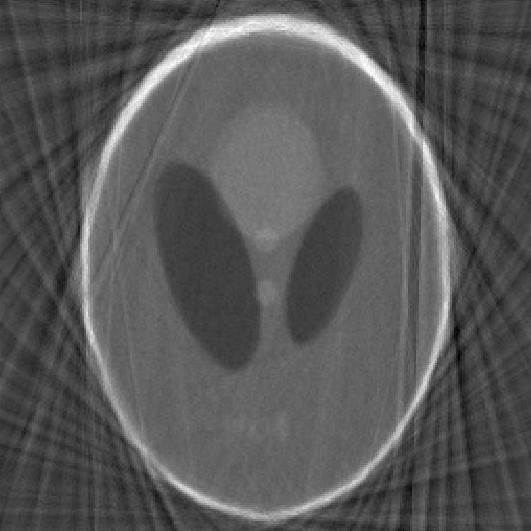}&
\includegraphics[width = 4cm]{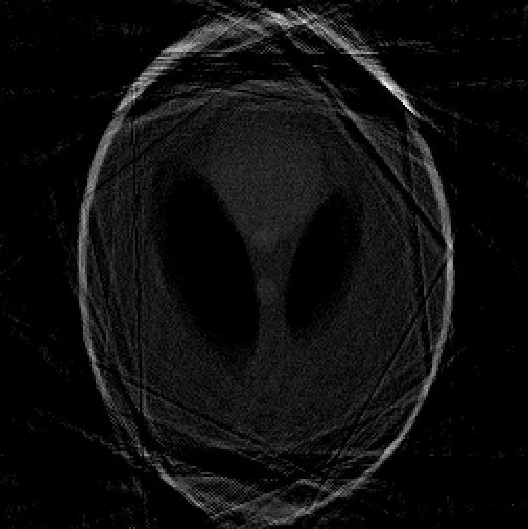}\\
{\footnotesize (a) True Solution} &
{\footnotesize (b) LSQR} &
{\footnotesize (c) IRN}&
{\footnotesize (d) FISTA}
\end{tabular}
\end{center}
\caption{Images for Linear Least Squares Comparison for the Second Test Problem} \label{fig:iterbraincomp2}
\end{figure}

Interestingly with a high perturbation level FISTA does the worst, while IRN continues to significantly outperform the LSQR algorithm. Unsurprisingly all three algorithms maintain a high relative error for the angles when larger perturbations are added. Figure \ref{fig:iterbraincomp2} shows the reconstructed image, and while no image is a high quality reconstruction, it is clear IRN is the cleanest. The two tests show though that there is certainly an image and perturbation level dependence as to which method has the best reconstructions.

\section{Conclusion and Outlook}\label{sec: conclusion}

We have devised an algorithm to effectively estimate the unknown geometry parameters of an uncalibrated portable CT machine. We exploited the problem structure to allow for parallelization and have shown significant speedup. Also, we have used fixed-point acceleration techniques to both reduce the image error and the number of iterations required for convergence. Additionally, we have surveyed and tested up-to-date methods to solve regularized linear least squares problems, and have demonstrated that the best choice of solver depends on the size of the perturbations to the geometry parameters.

In the future, we hope to investigate more about why the acceleration techniques work better only with the image vector. Additionally, our stopping criterion was chosen naively, and there may be better a method for choosing a stopping criterion and error tolerance, and we intend to investigate those. Finally, we will look into other acceleration techniques for the optimization problem. 

We also hope to apply our algorithm to other medical imaging applications where geometry parameters may be unknown such as bedside tomosynthesis \cite{CantEtal2017}.

\section{Acknowledgements}

We would like to express our appreciation to our mentor Dr. James Nagy. Additionally, we would like to thank the faculty at the Department of Mathematics, Emory University for dedicating their time to answering our questions, facilitating seminars and discussions related to this project. We also acknowledge the referees for their helpful comments in revising this manuscript. Lastly, we would like to acknowledge the National Science Foundation for giving us this opportunity to be part of NSF 2021 REU/RET Computational Mathematics for Data Science. This work was funded under NSF grant number DMS-2051019.

\bibliography{main}
\bibliographystyle{plain}

\begin{appendix}
\section{Running Numerical Experiments}\label{sec: Code Example}
Our implementation of algorithm \ref{alg: BCD} was built on top of the MATLAB IRtools package and naming conventions were chosen to mirror it. Our code can be found at https://github.com/manuelarturosantana/TomoREU2021. In the base IRtools package the function \texttt{PRset} is used to set up options for computed tomography problems as follows
$$\texttt{options = PRset(options, `field\_name1',field\_value1, `field\_name2',field\_value2)}.$$
We have updated \texttt{PRset} to accept values relating to BCD, such as an initial guess for the parameters. When a field name is not passed in default values are selected. Then to simulate a CT problem with unknown geometry parameters \texttt{PRtomo\_var} is used with the image size being $n \times n$
$$\texttt{[b,ProbInfo] = PRtomo\_var(n, options)}.$$

Above $\mathbf{b}$ is the right hand side vector with noise added, and \texttt{ProbInfo} contains the initial guess for the parameters. Thus, \texttt{PRtomo\_var} generates all the data necessary to simulate the inverse problem. Next the \texttt{IRset} function is used to set up hyper-parameters for the solvers in IRtools 

$$\texttt{options = IRset(options, `field\_name1',field\_value1, `field\_name2',field\_value2)}.$$
Again \texttt{IRset} was updated to include parameters for BCD, such as which acceleration technique to use. Finally a new function \texttt{IRbcd} computes the BCD
$$\texttt{[x,iterInfo] = IRbcd(b,iterOptions, probInfo)}.$$

\begin{figure}[H]
\begin{center}
\begin{tabular}{ccc}
\includegraphics[width = 5cm]{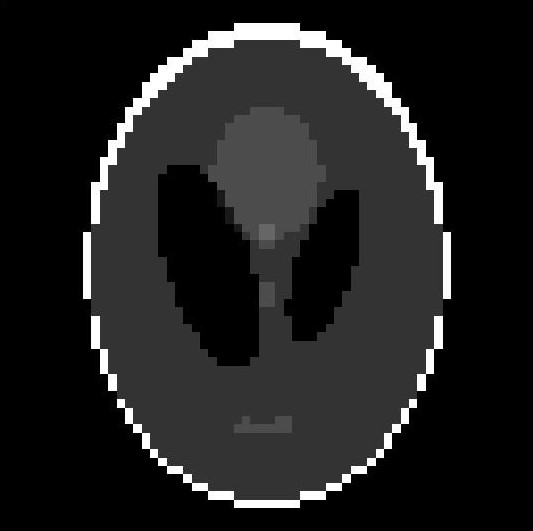}&
\includegraphics[width = 5cm]{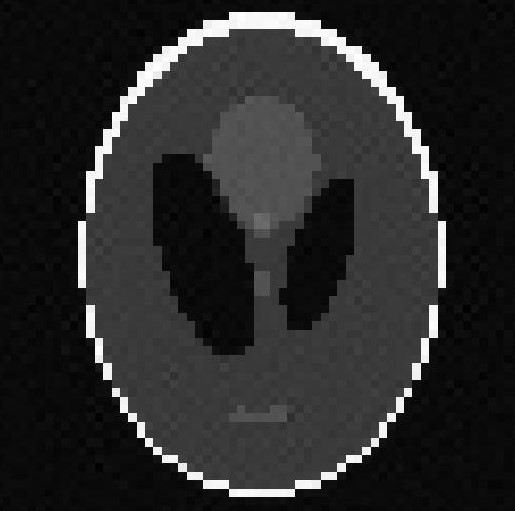}&
\includegraphics[width = 5cm]{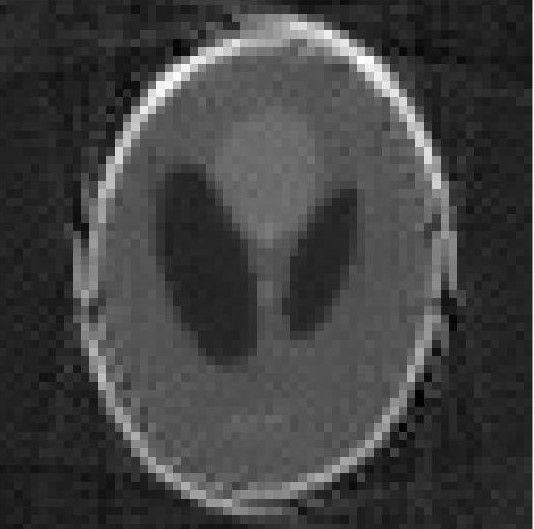}\\
{\footnotesize (a) True Solution} &
{\footnotesize (b) True Parameters Solution} &
{\footnotesize (c) Initial Parameters Solution}
\end{tabular}
\end{center}
\caption{Example Problem} \label{fig:Example problem}
\end{figure}

In the output \texttt{x} is the image vector after BCD terminates, and \texttt{iterInfo} contains information about errors at each iteration. Now we proceed to demonstrate an example of simulating and solving a CT reconstruction experiment with this code. To begin we generate the data
\begin{align*}
 &\texttt{n = 64;}\\
 &\texttt{ProbOptions = PRset(`Rpert',0.25,`anglespert',0.25);}\\
 &\texttt{[b,probInfo] = PRtomo\_var(n,ProbOptions);}
\end{align*}
This generates a test problem where perturbations are samplings of a uniform distribution on the interval $[-0.125,0.125]$ are added to the R values and the angles. For ease of computation the same perturbation is added to every fourth of the angles and R values. The image for this problem is the Shepp-Logan phantom. Figure \ref{fig:Example problem} $(a)$ represents the true solution, $(b)$ is the solution with noise in the sinogram $\mb{b}$, but true parameters, and $(c)$ is the solution with the initial guess parameters. 

The following codes sets up the iteration options and runs the problem 
\begin{align*}
    &\texttt{iterOptions= IRset(`nonlinSolver',`imfil',`accel',`anderson',`BCDmaxIter',10,...}\\
    &\texttt{`Rbounds',0.1250,`angleBounds',0.1250);}\\
    &\texttt{[x,iterInfo] = IRbcd(b,iterOptions,probInfo);}
\end{align*}

\begin{figure}[H]
\begin{center}
\begin{tabular}{cc}
\includegraphics[width = 5cm]{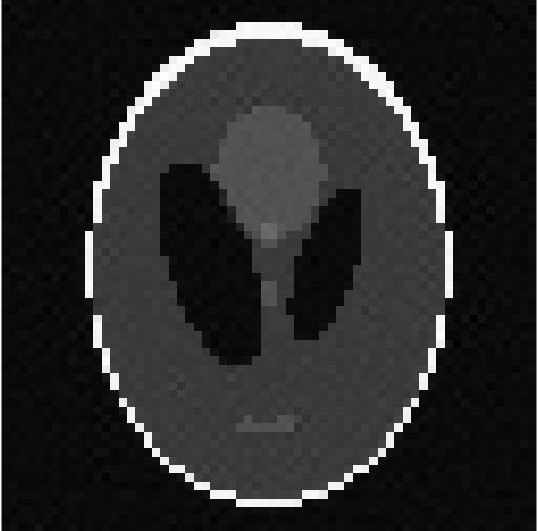}&
\includegraphics[width = 7cm]{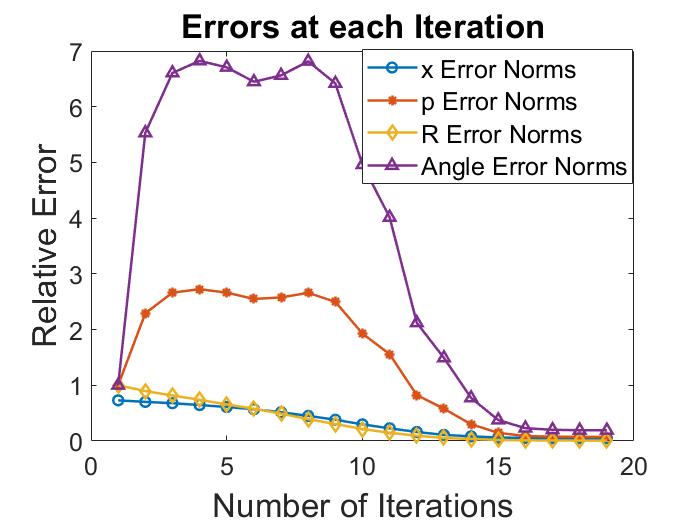}\\
{\footnotesize Image after BCD} &
\end{tabular}
\end{center}
\caption{Example Problem Solution} \label{fig:Example problem2}
\end{figure}

After solving the problem we see the solution and the relaltive error graph in figure \ref{fig:Example problem2}. This illustrates after the BCD algorithm we get a comparable solution to when the true geometry parameters are known.

The images in all five figures shown are generated with the following MATLAB code 
$$\texttt{PRshowbcd(iterInfo,probInfo)}.$$
\end{appendix}

\end{document}